\documentclass[12pt]{amsart}
\usepackage[cp850]{inputenc}
\usepackage{graphics,graphicx}

\newtheorem{theorem}{Theorem}[section]
\newtheorem{lemma}[theorem]{Lemma}
\newtheorem{proposition}[theorem]{Proposition}
\newtheorem{corollary}[theorem]{Corollary}

\theoremstyle{definition}

\newtheorem{example}[theorem]{Example}

\theoremstyle{remark}
\newtheorem{remark}[theorem]{Remark}

\usepackage{amssymb,amsmath}

\def\H{\mathbb{H}}

\newcommand{\m}{\mbox{$M$}}
\newcommand{\dm}{\mbox{$\mathrm{d}\Omega$}}
\newcommand{\s}{\mbox{$\Sigma$}}
\newcommand{\R}{\mbox{${\mathbb R}$}}
\newcommand{\N}{\mbox{$\m^2\times\R_1$}}
\newcommand{\g}[2]{\mbox{$\langle #1 ,#2 \rangle$}}
\newcommand{\fle}{\mbox{$\rightarrow$}}
\newcommand{\rf}[1]{\mbox{(\ref{#1})}}
\newcommand{\rl}[1]{{~\ref{#1}}}
\newcommand{\nablabar}{\mbox{$\overline{\nabla}$}}
\newcommand{\fs}{\mbox{$\mathcal{C}^\infty(\s)$}}
\newcommand{\f}{\mbox{$f:\Sigma^2\fle\N$}}

\def\beq{\begin{equation}}
\def\eeq{\end{equation}}

\begin{document}

\title[Calabi-Bernstein results for maximal surfaces]
{Calabi-Bernstein results for maximal surfaces in Lorentzian product spaces}

\author{Alma L. Albujer}
\address{Departamento de Matem\'{a}ticas, Universidad de Murcia, E-30100 Espinardo, Murcia, Spain}
\email{albujer@um.es}
\thanks{A.L. Albujer was supported by FPU Grant AP2004-4087 from Secretar\'\i a de Estado de Universidades e
Investigaci\'{o}n, MEC Spain.}

\author{Luis J. Al\'\i as}
\address{Departamento de Matem\'{a}ticas, Universidad de Murcia, E-30100 Espinardo, Murcia, Spain}
\email{ljalias@um.es}
\thanks{This work was partially supported by MEC project MTM2007-64504, and Fundaci\'{o}n S\'{e}neca
project 04540/GERM/06, Spain. This research is a result of the activity developed within the framework of the Programme in Support of Excellence Groups of the Regi\'{o}n de Murcia, Spain, by Fundaci\'{o}n S\'{e}neca, Regional Agency
for Science and Technology (Regional Plan for Science and Technology 2007-2010).}





\begin{abstract}
In this paper we establish new Calabi-Bernstein results for maximal surfaces immersed into a Lorentzian product space of
the form $\m^2\times\R_1$, where $\m^2$ is a connected Riemannian surface and $\m^2\times\R_1$ is endowed with the
Lorentzian metric $\g{}{}=\g{}{}_{M}-dt^2$. In particular, when \m\ is a Riemannian surface with
non-negative Gaussian curvature $K_M$, we prove that any complete maximal surface in \N\ must be
totally geodesic. Besides, if \m\ is non-flat we conclude that it must be a slice $\m\times\{t_0\}$, $t_0\in\R$
(here by \textit{complete} it is meant, as usual, that the induced Riemannian metric on the maximal surface from the ambient
Lorentzian metric is complete). We prove that the same happens if the maximal surface is complete with respect to the
metric induced from the Riemannian product $\m^2\times\R$. This allows us to give also a non-parametric version of the
Calabi-Bernstein theorem for entire maximal graphs in \N, under the same assumptions on $K_M$. Moreover,
we also construct counterexamples which show that our Calabi-Bernstein results are no longer true without
the hypothesis $K_M\geq 0$. These examples are constructed via a duality result between minimal and maximal graphs.
\end{abstract}

\maketitle

\section{Introduction}
\label{s1}
A maximal surface in a 3-dimensional Lorentzian manifold is a spacelike surface with zero mean curvature. Here by
\textit{spacelike} we mean that the induced metric from the ambient Lorentzian metric is a Riemannian metric on the
surface. The terminology \textit{maximal} comes from the fact that these surfaces locally maximize area among all
nearby surfaces having the same boundary \cite{Fr,BF}. Besides their mathematical interest, maximal surfaces and,
more generally, spacelike surfaces with constant mean curvature are also important in General Relativity (see, for
instance, \cite{MT}).

One of the most important global results about maximal surfaces is the Calabi-Bernstein theorem for maximal surfaces in
the 3-dimensional Lorentz-Minkowski space $\R^3_1$, which, in parametric version, states that the only complete maximal
surfaces in $\R^3_1$ are the spacelike planes. The Calabi-Bernstein theorem in $\R^3_1$ can be seen also in a
non-parametric form,
and it establishes that the only entire maximal graphs in $\R^3_1$ are the spacelike planes; that is, the only entire
solutions to the maximal surface equation
\[
\mathrm{Div}\left(\frac{Du}{\sqrt{1-|Du|^2}}\right)=0, \quad |Du|^2<1
\]
on the Euclidean plane $\R^2$ are affine functions.

This result was first proved by Calabi in \cite{Ca}, and extended to the general $n$-dimensional case by Cheng and Yau \cite{CY}. After that, several authors have approached the Calabi-Bernstein theorem for maximal surfaces ($n=2$) from different viewpoints, providing diverse extensions and new proofs of it, both in parametric and non-parametric versions
(Kobayashi \cite{Ko}, McNertey \cite{Mc}, Estudillo and Romero \cite{ER1,ER2,ER3}, Romero \cite{Ro}, Al\'\i as and
Palmer \cite{AP2}).

In this paper we establish new Calabi-Bernstein results for maximal surfaces immersed into a Lorentzian product space of the form $\m^2\times\R$, where $\m^2$ is a connected Riemannian surface and $\m^2\times\R$ is endowed with the Lorentzian metric
$$
\g{}{}=\pi_{M}^*(\g{}{}_{M})-\pi_{\mathbb{R}}^*(dt^2).
$$
Here $\pi_{M}$ and $\pi_{\mathbb{R}}$ denote the projections from $\m\times\R$ onto each factor, and $\g{}{}_{M}$ is the
Riemannian metric on \m. For simplicity, we will write simply
$$
\g{}{}=\g{}{}_{M}-dt^2,
$$
and we will denote by \N\ the 3-dimensional product manifold $\m^2\times\R$ endowed with that Lorentzian metric. In particular, when $\m^2=\R^2$ is the flat Euclidean plane, then $\N=\R^3_1$ is nothing but the Lorentz-Minkowski space.
Before describing our results, it is worth pointing out that the Calabi-Bernstein result is not true in general Lorentzian product spaces; counterexamples to it have been given recently by the first author in \cite{Al}, when $\m=\H^2$
is the hyperbolic plane. In Section\rl{dualitysection} we also construct new counterexamples (see below). Other
examples of non-trivial complete maximal surfaces in $\H^2\times\R_1$ are given in \cite{FM}.

Our first main result is Theorem\rl{CBth}, which states that any complete maximal surface \s\ immersed into a Lorentzian
product \N, where \m\ is a (necessarily complete) Riemannian surface with non-negative Gaussian curvature, must be
totally geodesic. Moreover, if \m\ is non-flat we conclude that \s\ must be a slice $\m\times\{t_0\}$, $t_0\in\R$. Here
by \textit{complete} it is meant, as usual, that the induced Riemannian metric on \s\ from the ambient Lorentzian metric
is complete. In Theorem\rl{CBth3} we prove that the same happens if \s\ is complete with respect to the metric induced
from the Riemannian product $\m^2\times\R$. This allows us to give a non-parametric version of the Calabi-Bernstein
theorem (Theorem\rl{CBth2} and Corollary\rl{CBco2}), where we prove that any entire maximal graph in \N\ must be totally
geodesic and conclude that the only entire solutions to the maximal surface equation on any complete non-flat Riemannian
surface \m\ with non-negative Gaussian curvature are the constant functions.

As observed above, in all of these results, the assumption on the Gaussian curvature of \m\ is necessary as shown by the fact that when
$\m=\H^2$ is the hyperbolic plane, then there exist examples of complete maximal surfaces in $\H^2\times\R_1$ which are
not totally geodesic, as well as examples of non-trivial entire maximal graphs over $\H^2$, see Examples\rl{ex1}
and\rl{ex2} (other examples have been given recently by the first author in \cite{Al}, where
explicit examples of non-trivial entire maximal graphs in $\H^2\times\R_1$ have been found by looking for explicit solutions
of the corresponding partial differential equation on $\H^2$). Particularly interesting is the fact that the entire maximal
graph given in Example\rl{ex2} is not complete. As is well known, such circumstance cannot occur in the
Lorentz-Minkowski space $\mathbb{R}^3_1$, since by a result of Cheng and Yau \cite{CY}, surfaces with constant mean curvature and
closed in $\mathbb{R}^3_1$ are necessarily complete. For the construction of our Examples \rl{ex1} and\rl{ex2}, in Theorem\rl{duality} we establish a
simple but nice duality result between solutions to the minimal surface equation in a Riemannian product $\m^2\times\R$
and solutions to the maximal surface equation in a Lorentzian product $\m^2\times\R_1$. This extends
\cite[Theorem 3]{AP1}.

\section{Preliminaries}
A smooth immersion \f\ of a connected surface $\Sigma^2$ is said to be a spacelike surface if $f$ induces a Riemannian
metric on $\Sigma$, which as usual is also denoted by $\g{}{}$. In that case, since
$$
\partial_t=(\partial/\partial_t)_{(x,t)}, \quad x\in\m, t\in\R,
$$
is a unitary timelike vector field globally defined on the ambient spacetime \N, then there exists a unique unitary
timelike normal field $N$ globally defined on $\Sigma$ which is in the same time-orientation as $\partial_t$, so that
$$
\g{N}{\partial_t}\leq -1<0 \quad \mathrm{on} \quad \Sigma.
$$
We will refer to $N$ as the future-pointing Gauss map of $\Sigma$, and we will denote by
$\Theta:\s\fle (-\infty,-1]$ the smooth function on \s\ given by $\Theta=\g{N}{\partial_t}$. Observe that the function
$\Theta$ measures the hyperbolic angle $\theta$ between the future-pointing vector fields $N$ and $\partial_t$ along \s. Indeed, they are related by $\cosh{\theta}=-\Theta$.

Let $\nablabar$ and $\nabla$ denote the Levi-Civita connections in \N\ and $\Sigma$, respectively. Then the Gauss and
Weingarten formulae for the spacelike surface \f\ are given by
\beq
\label{gaussfor}
\nablabar_XY=\nabla_XY-\g{AX}{Y}N
\eeq
and
\beq
\label{wein}
AX=-\nablabar_XN,
\eeq
for any tangent vector fields $X,Y\in T\Sigma$. Here $A:T\Sigma\fle T\Sigma$ stands for the shape operator (or second fundamental form) of $\Sigma$ with respect to its future-pointing Gauss map $N$. As is well known, the Gaussian curvature
$K$ of the surface \s\ is described in terms of $A$ and the curvature of the ambient spacetime by the Gauss
equation, which is given by
\beq
\label{gausseq}
K=\overline{K}-\mathrm{det}A,
\eeq
where $\overline{K}$ denotes the sectional curvature in \N\ of the plane tangent to \s. On the other hand, it is not
difficult to see that the curvature tensor $\overline{R}$ of \N\ can be written in terms of the Gaussian curvature of \m\
by
\begin{eqnarray}
\label{curvature}
\nonumber \overline{R}(U,V)W & = & K_M(\pi_M)\left(\g{U}{W}V-\g{V}{W}U\right)+\\
{} & {} & K_M(\pi_M)\g{W}{\partial_t}\left(\g{U}{\partial_t}V-\g{V}{\partial_t}U\right)+\\
\nonumber {} & {} & K_M(\pi_M)\left(\g{U}{W}\g{V}{\partial_t}-\g{V}{W}\g{U}{\partial_t}\right)\partial_t.
\end{eqnarray}
for any vector fields $U,V,W$ tangent to \N.
In particular, if $\{ E_1,E_2\}$ is a local orthonormal frame on \s, then
from \rf{curvature} we obtain that
\beq
\label{eq1}
\overline{K}=\g{\overline{R}(E_1,E_2)E_1}{E_2}=\kappa_M(1+\g{E_1}{\partial_t^\top}^2+\g{E_2}{\partial_t^\top}^2)=
\kappa_M(1+\|\partial_t^\top\|^2),
\eeq
where, for simplicity, $\kappa_M$ stands for the Gaussian curvature of \m\ along the surface \s, that is,
$\kappa_M=K_M\circ\Pi\in\fs$ where $\Pi=\pi_M\circ f:\Sigma\fle M$ denotes the projection of \s\ onto \m. Here and
in what follows, if $Z$ is a vector field along the immersion \f, then $Z^\top\in T\Sigma$ denotes the tangential
component of $Z$ along \s, that is,
\[
Z=Z^\top-\g{N}{Z} N,
\]
and $\|\cdot\|$ denotes the norm of a vector field on \s.
In particular, $\partial_t^\top=\partial_t+\Theta N$ and then
\beq
\label{eq2}
-1=\|\partial_t^\top\|^2-\Theta^2.
\eeq
Therefore, expression \rf{eq1} becomes $\overline{K}=\kappa_M\Theta^2$, and the Gauss equation \rf{gausseq} can be written as
\beq
\label{gaussEq}
K=\kappa_M\Theta^2-\mathrm{det}A.
\eeq
On the other hand, the Codazzi equation of the spacelike surface \s\ describes the
tangent component of $\overline{R}(X,Y)N$, for any tangent vector fields $X,Y\in T\Sigma$, in terms of the derivative of
the shape operator, and it is given by
\beq
\label{codazzi}
(\overline{R}(X,Y)N)^\top=(\nabla_XA)Y-(\nabla_YA)X.
\eeq
Here $\nabla_XA$ denotes the covariant derivative of $A$, that is,
$$
(\nabla_XA)Y=\nabla_X(AY)-A(\nabla_XY).
$$
By \rf{curvature}, Codazzi equation \rf{codazzi} becomes
\beq
\label{codazzibis}
(\nabla_XA)Y=(\nabla_YA)X+\kappa_M\Theta\left(\g{X}{\partial_t^\top}Y-\g{Y}{\partial_t^\top}X\right).
\eeq

For a spacelike surface \f, we will call the height function of \s, denoted by $h$, the projection of \s\ onto
\R, that is, $h\in\fs$ is the smooth function given by $h=\pi_{\mathbb{R}}\circ f$. Observe that the gradient of $\pi_{\mathbb{R}}$ on \N\ is
\[
\nablabar\pi_{\mathbb{R}}=-\g{\nablabar\pi_{\mathbb{R}}}{\partial_t}\partial_t=-\partial_t.
\]
Therefore, the gradient of $h$ on \s\ is
\[
\nabla h=(\nablabar\pi_{\mathbb{R}})^\top=-\partial_t^\top.
\]
Observe that from \rf{eq2} we have
\beq
\label{eq6}
\|\nabla h\|^2=\Theta^2-1.
\eeq
Since $\partial_t$ is parallel on \N\ we have that
\beq
\label{eq4}
\nablabar_X\partial_t=0
\eeq
for any tangent vector field $X\in T\Sigma$. Writing $\partial_t=-\nabla h-\Theta N$ along the surface \s\ and using
Gauss \rf{gaussfor} and Weingarten \rf{wein} formulae, we easily get from \rf{eq4} that
\beq
\label{eq5}
\nabla_X\nabla h=\Theta AX
\eeq
for every $X\in T\Sigma$. Therefore the Laplacian on \s\ of the height function is given by
\beq
\label{laplah}
\Delta h=\Theta\mathrm{tr}A=-2H\Theta,
\eeq
where $H=-(1/2)\mathrm{tr}A$ is the mean curvature of \s\ relative to $N$.

On the other hand, \rf{eq4} also yields
\[
X(\Theta)=\g{\nablabar_XN}{\partial_t}=\g{AX}{\nabla h}=\g{X}{A\nabla h}
\]
for every $X\in T\Sigma$, and then the gradient of $\Theta$ on \s\ is given by
\beq
\label{gradT}
\nabla\Theta=A\nabla h.
\eeq
From here, using Codazzi equation \rf{codazzibis} and equations \rf{eq6} and \rf{eq5} we get
\begin{eqnarray*}
\nabla_X\nabla\Theta & = & (\nabla_XA)(\nabla h)+A(\nabla_X\nabla h) \\
{} & = & (\nabla_{\nabla h}A)(X)+\kappa_M\Theta\left(\g{X}{\partial_t^\top}\nabla h-\g{\nabla h}{\partial_t^\top}X\right)
+\Theta A^2X \\
{} & = & (\nabla_{\nabla h}A)(X)+\kappa_M\Theta\left((\Theta^2-1)X-\g{X}{\nabla h}\nabla h\right)
+\Theta A^2X,
\end{eqnarray*}
for every $X\in T\Sigma$. Thus, the Laplacian of $\Theta$ is given by
\begin{eqnarray}
\label{laplaT}
\Delta\Theta & = & \mathrm{tr}(\nabla_{\nabla h}A)+\kappa_M\Theta(\Theta^2-1)+\Theta\mathrm{tr}(A^2) \\
\nonumber {} & = & -2\g{\nabla H}{\nabla h}+\Theta(\kappa_M(\Theta^2-1)+\|A\|^2),
\end{eqnarray}
with $\|A\|^2=\mathrm{tr}(A^2)$, where we are using the fact that
\[
\mathrm{tr}(\nabla_{\nabla h}A)=\nabla_{\nabla h}(\mathrm{tr}A)=-2\g{\nabla H}{\nabla h}.
\]

A spacelike surface \s\ is said to be a maximal surface if its mean curvature vanishes, $H=0$ on \s. Equation \rf{laplah}
implies that if \s\ is a maximal surface in \N, then $h$ is a harmonic function, $\Delta h=0$ on \s.
Besides, if \s\ is maximal then
\beq
\label{eq7}
A^2=\frac{1}{2}\|A\|^2I,
\eeq
where $I$ denotes the identity map on $T\Sigma$ and $\|A\|^2=-2\mathrm{det}A$. Therefore, the Gauss equation \rf{gaussEq}
becomes
\beq
\label{gaussmaximal}
K=\kappa_M\Theta^2+\frac{1}{2}\|A\|^2.
\eeq
On the other hand, from \rf{eq6}, \rf{gradT} and \rf{eq7} we also obtain that for a maximal surface it holds that
\beq
\label{normamaximal}
\|\nabla\Theta\|^2=\frac{1}{2}\|A\|^2\|\nabla h\|^2=\frac{1}{2}\|A\|^2(\Theta^2-1).
\eeq

\section{A Calabi-Bernstein theorem for maximal surfaces}
\label{mecagoenschoen}
We start by stating the following remarkable property.
\begin{lemma}
\label{lemmaonto}
Let $\m^2$ be a Riemannian surface. If \N\ admits a complete spacelike surface \f, then \m\ is necessarily complete and
the projection $\Pi=\pi_M\circ f:\Sigma\fle M$ is a covering map.
\end{lemma}
\begin{proof}
We follow the ideas in the proof of \cite[Lemma 3.1]{ARS1}.
Let \f\ be a spacelike surface and consider $\Pi=\pi_M\circ f:\Sigma\fle M$ its projection on \m. It is not difficult
to see that $\Pi^*(\g{}{}_{M})\geq\g{}{}$, where $\g{}{}$ stands for the Riemannian metric on \s\ induced
from the Lorentzian ambient space. That means that $\Pi$ is a local diffeomorphism which increases the distance
between the Riemannian surfaces \s\ and \m. Then, the proof finishes recalling that if a map, from a connected complete
Riemannian manifold $M_1$ into another connected Riemannian manifold $M_2$ of the same dimension, increases the distance,
then it is a covering map and $M_2$ is complete \cite[Chapter VIII, Lemma 8.1]{KN}.
\end{proof}
In particular, if \N\ admits a compact spacelike surface, then \m\ is necessarily compact (see
\cite[Proposition 3.2 (i)]{ARS1}). An easy consequence of \rf{laplah} is the following.
\begin{proposition}
\label{prop3.2}
Let $\m^2$ be a Riemannian surface. If \f\ is a compact spacelike surface in \N\ whose mean
curvature $H$ does not change sign, then it must be a slice $\m\times\{ t_0\}$, $t_0\in\R$.
In particular, the only compact maximal surfaces in \N\ are the slices.
\end{proposition}
For the proof simply observe that, since $\Theta<0$ and $H$ does not change sign, then formula \rf{laplah} says
that the height function $h$ must be either subharmonic or superharmonic on $\Sigma$, according to the sign of $H$. But the compactness
of $\Sigma$ implies that $h$ must be constant, and the conclusion of Proposition\rl{prop3.2} follows.

Under completeness assumption, we have the following parametric version of a Calabi-Bernstein result in \N.
\begin{theorem}
\label{CBth}
Let $\m^2$ be a (necessarily complete) Riemannian surface with non-negative Gaussian curvature, $K_M\geq 0$. Then any
complete maximal surface $\Sigma^2$ in \N\ is totally geodesic. In addition, if $K_M>0$ at some point on \m, then \s\ is
a slice $\m\times\{ t_0\}$, $t_0\in\R$.
\end{theorem}
As a direct consequence of Theorem\rl{CBth} we have the following.
\begin{corollary}
\label{CBco}
Let $\m^2$ be a complete non-flat Riemannian surface with non-negative Gaussian curvature, $K_M\geq 0$. Then the
only complete maximal surfaces in \N\ are the slices $\m\times\{ t_0\}$, $t_0\in\R$.
\end{corollary}

Observe that if $\m^2=\R^2$ is the flat Euclidean plane, then $\N=\R^3_1$ is nothing but the 3-dimensional
Lorentz-Minkowski space, and any spacelike affine plane in $\R^3_1$ which is not horizontal determines a complete
totally geodesic surface which is not a slice. On the other hand, the assumption $K_M\geq 0$ is necessary as shown by
the fact that there exist examples of non-totally geodesic complete maximal surfaces in $\mathbb{H}^2\times\R_1$,
where $\mathbb{H}^2$ is the hyperbolic plane (see Example\rl{ex1}).

\begin{proof}[Proof of Theorem\rl{CBth}]
Since \s\ is maximal and $\kappa_M\geq 0$, \rf{gaussmaximal} implies that $K\geq 0$ on \s. Then, \s\ is a complete Riemannian surface
with non-negative Gaussian curvature and, by a classical result due to Ahlfors \cite{Ah} and
Blanc-Fiala-Huber \cite{Hu}, we know that \s\ is parabolic, in the sense that any non-positive subharmonic function on
the surface must be constant.

Recall that $\Theta\leq -1<0$. From \rf{laplaT} and \rf{normamaximal}, we can compute
\beq
\label{laplaT2}
\Delta\left(\frac{1}{\Theta}\right)=-\frac{\Delta\Theta}{\Theta^2}+\frac{2\|\nabla\Theta\|^2}{\Theta^3}=
-\frac{1}{\Theta}\left(\kappa_M(\Theta^2-1)+\frac{\|A\|^2}{\Theta^2}\right)\geq 0.
\eeq
That is, $1/\Theta$ is a negative subharmonic function on the parabolic surface \s, and hence it must be constant. That is,
$\Theta=\Theta_0\leq -1$ is constant on \s, and by \rf{laplaT2} we also get that $\|A\|^2=0$ and
$\kappa_M(\Theta_0^2-1)=0$ on \s. Therefore, \s\ is totally geodesic in \N\ and, if $\kappa_M>0$ at some point on \s,
then it must be $\Theta_0=-1$, which by \rf{eq6} means that $h$ is constant and \s\ is a slice. Finally, observe that
since the projection $\Pi:\Sigma\fle M$ is onto (Lemma\rl{lemmaonto}), then $\kappa_M>0$ at some point on \s\ if and only
if $K_M>0$ at some point on \m.
\end{proof}

\section{Entire maximal graphs and Calabi-Bernstein theorem}
Let $\Omega\subseteq\m^2$ be a connected domain. Every smooth function $u\in\mathcal{C}^\infty(\Omega)$ determines a
graph over $\Omega$ given by $\Sigma(u)=\{ (x,u(x)) : x\in\Omega \}\subset\N$. The metric induced on $\Omega$ from the
Lorentzian metric on the ambient space via $\Sigma(u)$ is given by
\beq
\label{gu}
\g{}{}=\g{}{}_M-du^2.
\eeq
Therefore, $\Sigma(u)$ is a spacelike surface in \N\ if and only if $|Du|^2<1$ everywhere on $\Omega$, where $Du$ denotes the gradient of $u$ in $\Omega$ and $|Du|$ denotes its norm, both with respect to the original metric $\g{}{}_M$ on $\Omega$. If $\Sigma(u)$ is a spacelike graph over a domain $\Omega$, then it is not difficult to see that the vector field
\beq
\label{normal}
N(x)=\frac{1}{\sqrt{1-|Du(x)|^2}}\left(Du(x)+\partial_t|_{(x,u(x))}\right), \quad x\in\Omega,
\eeq
defines the future-pointing Gauss map of $\Sigma(u)$. The shape operator of $\Sigma(u)$ with respect to $N$ is given by
\beq
\label{sff}
AX=-\frac{1}{\sqrt{1-|Du|^2}}D_XDu-\frac{\g{D_XDu}{Du}_M}{(1-|Du|^2)^{3/2}}Du,
\eeq
for every tangent vector field $X$ on $\Omega$, where $D$ denotes the Levi-Civita connection in $\Omega$ with respect to the metric $\g{}{}_M$. It follows from here that the mean curvature $H(u)$ of a spacelike graph $\Sigma(u)$ is given by
\[
2H(u)=\mathrm{Div}\left(\frac{Du}{\sqrt{1-|Du|^2}}\right)
\]
where Div stands for the divergence operator on $\Omega$ with respect to the metric $\g{}{}_M$. In particular, $\Sigma(u)$ is a maximal graph if and only if the function $u$ satisfies the following partial differential equation on the domain $\Omega$,
\beq
\label{zmc}
\mathrm{Maximal}[u]=\mathrm{Div}\left(\frac{Du}{\sqrt{1-|Du|^2}}\right)=0, \quad |Du|^2<1.
\eeq

A graph is said to be entire if $\Omega=\m$.
As a direct consequence of Lemma\rl{lemmaonto}, it follows that when \m\ is a complete Riemannian surface which is simply connected, then every complete spacelike surface in \N\ is an entire graph. In fact, since \m\ is
simply connected then the projection $\Pi$ is a diffeomorphism between \s\ and \m, and hence \s\ can be written as the
graph over \m\ of the function $u=h\circ\Pi^{-1}\in\mathcal{C}^\infty(\m)$. However, in contrast to the case of graphs
into a Riemannian product space, an entire spacelike graph in \N\ is not necessarily complete, in the sense that the
induced Riemannian metric \rf{gu} is not necessarily complete on \m. For instance, let $u:\R^2\fle\R$ be a function defined by
\[
u(x_1,x_2)=\int_0^{|x_1|}\sqrt{1-e^{-s}}ds
\]
when $|x_1|\geq 1$, and $u(x_1,x_2)=\phi(x_1)$ when $|x_1|<1$, where $\phi\in\mathcal{C}^\infty(\R)$ is a smooth
extension satisfying $\phi'(s)^2<1$ for all $s\in(-1,1)$. Then $u$ determines an entire spacelike graph $\Sigma(u)$ in
$\R^2\times\R_1=\R^3_1$ which is not complete. In fact, observe that the curve $\alpha:\R\fle\Sigma(u)$ given by
$\alpha(s)=(s,0,u(s,0))$ is a divergent curve in $\Sigma(u)$ with finite length, because of
\[
\int_{-\infty}^{+\infty}\|\alpha'(s)\|ds=
\int_{-1}^{1}\sqrt{1-\phi'(s)^2}ds+2\int_1^{+\infty}e^{-s/2}ds<2(1+2/\sqrt{e}).
\]
As another example, in Example\rl{ex2} we construct an example of an entire maximal graph in $\H^2\times\R_1$ which is
not complete.

For that reason, the Calabi-Bernstein result given at Theorem\rl{CBth} does not imply in principle that, under the same
hypothesis on \m, any entire maximal graph in \N\ must be totally geodesic. This is certainly true for entire maximal
graphs in the Lorentz-Minkowski space \cite{Ro}, and also for entire maximal graphs in Robertson-Walker spaces of the
form $\R^2\times_\varrho\R_1$, under certain assumptions on the warping function $\varrho$ (for the details, see
\cite{LR1}). However, although we cannot establish a similar result for entire maximal graphs in \N\ just as a direct
consequence of our Theorem\rl{CBth}, we can obtain it as a consequence of the following result.
\begin{theorem}
\label{CBth3}
Let $\m^2$ be a (non necessarily complete) Riemannian surface with non-negative Gaussian curvature, $K_M\geq 0$. Then any
maximal surface $\Sigma^2$ in \N\ which is complete with respect to the metric induced from the Riemannian product
$\m^2\times\R$ is totally geodesic. In addition, if $K_M>0$ at some point on $\Sigma$, then \m\ is necessarily complete and
$\Sigma$ is a slice $\m\times\{t_0\}$.
\end{theorem}
In particular, if $\m^2$ is complete and \f\ is a spacelike surface which is properly immersed in \N, then the
metric induced on \s\ from the Riemannian product $\m^2\times\R$ is complete. Then, we have the following consequence.
\begin{corollary}
\label{CBco3}
Let $\m^2$ be a complete Riemannian surface with non-negative Gaussian curvature, $K_M\geq 0$. Then any
maximal surface $\Sigma^2$ properly immersed into \N\ is totally geodesic. In addition, if $K_M>0$ at some point on
$\Sigma$, then $\Sigma$ is a slice $\m\times\{t_0\}$.
\end{corollary}
This happens, for instance, when $\m^2$ is complete and $\s\subset\N$ is a closed embedded maximal surface. In particular,
it happens for entire maximal graphs, and it yields the following non-parametric version of the Calabi-Bernstein theorem.
\begin{theorem}
\label{CBth2}
Let $\m^2$ be a complete Riemannian surface with non-negative Gaussian curvature, $K_M\geq 0$. Then any
entire maximal graph $\Sigma(u)$ in \N\ is totally geodesic. In addition, if $K_M>0$ at some point on \m, then $u$ is
constant.
\end{theorem}
As a consequence of this we also have the following.
\begin{corollary}
\label{CBco2}
Let $\m^2$ be a complete non-flat Riemannian surface with non-negative Gaussian curvature, $K_M\geq 0$. Then the
only entire solutions to the maximal surface equation \rf{zmc} are the constant functions.
\end{corollary}
The proof of Theorem\rl{CBth3} follows the ideas introduced by Romero in his proof of \cite[Theorem]{Ro}
(see also the proof of \cite[Theorem A]{LR1}), which was inspired by Chern's proof of the classical Bernstein theorem for
entire minimal graphs in Euclidean space \cite{Ch}. Examples\rl{ex1} and\rl{ex2} show that the assumption $K_M\geq 0$ is
necessary.

\begin{proof}[Proof of Theorem\rl{CBth3}]
Let \s\ be a maximal surface in \N. For simplicity, we denote by $g=\g{}{}$ the Riemannian metric induced on \s\ from the Lorentzian product \N. Since $1-\Theta\geq 2>0$, we may introduce on \s\ the conformal metric
\[
\hat{g}=(1-\Theta)^2g.
\]
As is well known, the Gaussian curvature $\hat{K}$ of $(\s,\hat{g})$ is given by
\beq
\label{khat}
(1-\Theta)^2\hat{K}=K-\Delta\log{(1-\Theta)},
\eeq
where $K$ is the Gaussian curvature of $(\s,g)$, which is given by \rf{gaussmaximal}. Using \rf{laplaT} and \rf{normamaximal}, we can compute
\[
\Delta\log{(1-\Theta)}=-\frac{\Delta\Theta}{1-\Theta}-\frac{\|\nabla\Theta\|^2}{(1-\Theta)^2}=
\frac{1}{2}\|A\|^2+\Theta(\Theta+1)\kappa_M,
\]
which by \rf{gaussmaximal} becomes
\beq
\label{laplalog}
\Delta\log{(1-\Theta)}=K+\Theta\kappa_M.
\eeq
Therefore, from \rf{khat} we conclude that $\hat{K}\geq 0$ on \s. On the other hand, we also have that
\beq
\label{eq8}
\hat{g}(X,X)=(1-\Theta)^2g(X,X)\geq\Theta^2g(X,X)=\Theta^2|X^*|^2-\Theta^2X(h)^2
\eeq
for every tangent vector field $X$ on \s. Here $U^*=(\pi_{M})_*(U)$ denotes projection
onto the surface $M^2$ of a vector field $U$ defined on \N, that is,
\[
U=U^*-\g{U}{\partial_t}\partial_t,
\]
and we recall that $|\cdot |$ denotes the norm with respect to the original metric $\g{}{}_M$ on \m. Writing
\beq
\label{eq8.5}
X=X^*-\g{X}{\partial_t}\partial_t=X^*+X(h)\partial_t
\eeq
and
\[
N=N^*-\g{N}{\partial_t}\partial_t=N^*-\Theta\partial_t,
\]
we find that
\[
|N^*|^2=\Theta^2-1 \quad \mbox{and} \quad \Theta X(h)=-\g{X^*}{N^*}_M.
\]
Then, by Cauchy-Schwarz inequality we have
\[
\Theta^2X(h)^2=\g{X^*}{N^*}_M^2\leq |X^*|^2|N^*|^2=|X^*|^2(\Theta^2-1),
\]
which jointly with \rf{eq8} yields
\beq
\label{eq8.55}
\hat{g}(X,X)\geq |X^*|^2=\g{X^*}{X^*}_M.
\eeq
Let $g'$ denote the Riemannian metric induced on \s\ from the Riemannian product $\m^2\times\R$. From \rf{eq8.5} we have
\[
g(X,X)=|X^*|^2-X(h)^2 \quad \mbox{and} \quad g'(X,X)=|X^*|^2+X(h)^2.
\]
Therefore,
\[
|X^*|^2=\frac{1}{2}(g(X,X)+g'(X,X))\geq\frac{1}{2}g'(X,X),
\]
and by \rf{eq8.55} we get
\[
\hat{g}(X,X)\geq |X^*|^2\geq\frac{1}{2}g'(X,X)
\]
for every tangent vector field $X$ on \s. This implies that $\hat{L}\geq (1/\sqrt{2})L'$, where $\hat{L}$ and $L'$ denote the length of a curve on \s\ with respect to the Riemannian metrics $\hat{g}$ and $g'$, respectively. As a consequence, since
we are assuming that the metric $g'$ is complete on \s, then $\hat{g}$ is also complete.

Summing up, $(\s^2,\hat{g})$ is a complete Riemannian surface with non-negative Gaussian curvature and, from the same classical result by Ahlfors and Blanc-Fiala-Huber used in the proof of Theorem\rl{CBth}, we conclude that $(\s^2,\hat{g})$ is parabolic. Since $n=2$, the Laplacian $\Delta$ on \s\ with respect to $g$ and the Laplacian $\hat{\Delta}$ on \s\ with respect to the conformal metric $\hat{g}$ are related by
\[
\Delta=(1-\Theta)^2\hat{\Delta},
\]
which implies that the property of being subharmonic is preserved under conformal changes of metric. Therefore,
$(\s^2,g)$ is also parabolic. The proof then follows as in the proof of Theorem\rl{CBth}, since by \rf{laplaT2} we know that $1/\Theta$ is a negative subharmonic function on $(\s^2,g)$.
\end{proof}

\begin{remark}
It is worth pointing out that Theorem\rl{CBth} can be also seen as a consequence of Theorem\rl{CBth3}, because
every complete spacelike surface $\Sigma$ in \N\ is also complete with respect to the metric induced from the Riemannian product
$M^2\times\mathbb{R}$. This follows from the fact that $g'\geq g$, where $g$ and $g'$ stand for the metrics induced
on $\Sigma$ from the Lorentzian product and the Riemannian product, respectively. Nevertheless, the proof of
Theorem\rl{CBth} given in Section\rl{mecagoenschoen} is much simpler and direct.
\end{remark}

\section{A duality result between minimal and maximal graphs}
\label{dualitysection}
In \cite{Ca} Calabi observed a simple but nice duality between solutions to the
minimal surface equation in the Euclidean space $\R^3$ and solutions to the maximal surface equation in the
Lorentz-Minkowski space $\R^3_1$ (see also \cite{AP1} for an alternative approach to that duality given by the second
author, jointly with Palmer). By regarding $\R^3$ as the Riemannian product space $\R^2\times\R$ and $\R^3_1$ as the
Lorentzian product space $\R^2\times\R_1$, we observe here that the same duality holds in general between solutions to
the minimal surface equation in a Riemannian product space $\m\times\R$ and solutions to the maximal surface equation in
a Lorentzian product space $\m\times\R_1$. First of all, recall that a smooth function $u$ on a connected domain
$\Omega\subseteq\m^2$ defines a minimal graph $\Sigma(u)$ in $\m\times\R$ if and only if $u$ satisfies the following
partial differential equation on $\Omega$,
\beq
\label{minimal}
\mathrm{Minimal}[u]=\mathrm{Div}\left(\frac{Du}{\sqrt{1+|Du|^2}}\right)=0,
\eeq
where, as in \rf{zmc}, Div and $Du$ stand for the divergence operator and the gradient of $u$ in $\Omega$ with respect to the metric $\g{}{}_M$, respectively. Following the ideas in \cite{AP1} we can prove the following general result.
\begin{theorem}
\label{duality}
Let $\Omega\subseteq\m^2$ be a simply connected domain of a Riemannian surface $\m^2$. There exists a non-trivial solution $u$ to the minimal surface equation on $\Omega$,
\[
\mathrm{Minimal}[u]=0,
\]
if and only if there exists a non-trivial solution $w$ to the maximal surface equation on $\Omega$,
\[
\mathrm{Maximal}[w]=0, \quad |Dw|^2<1.
\]
\end{theorem}
Here by a \textit{non-trivial solution} we mean a solution with non-parallel gradient. Observe that, from \rf{sff}, a
spacelike graph determined by a smooth function $u\in\mathcal{C}^\infty(\Omega)$ is totally geodesic in \N\ if and only
if $Du$ is parallel on $\Omega$. Similarly, a graph determined by a function $u$ is totally geodesic in $\m^2\times\R$ if
and only if $Du$ is parallel. Therefore, non-trivial solutions to either the minimal or maximal surface equation
correspond to non-totally geodesic either minimal or maximal graphs.

\begin{proof}
Since $\Omega$ is simply connected, it is orientable and can be endowed with a globally defined area form \dm\
and an almost complex structure $J$. Recall that, for every vector field $X$ on $\Omega$, it holds that
\beq
\label{eq9}
\mathrm{Div}X\,\dm=d\omega_{JX}
\eeq
where $\omega_{JX}$ denotes the 1-form in $\Omega$ which is metrically equivalent to the field $JX$, that is,
\[
\omega_{JX}(Y)=\g{JX}{Y}_M.
\]
Now the proof of Theorem\rl{duality} follows as the proof of \cite[Theorem 3]{AP1}. For the sake of completeness we sketch it here. Assume that $u$ is a non-trivial solution of \rf{minimal} on the domain $\Omega$. Then, by \rf{eq9} the 1-form $\omega_{JU}$ is closed on $\Omega$, where $U$ is the vector field on $\Omega$ given by
\[
U=\frac{Du}{\sqrt{1+|Du|^2}}.
\]
Since $\Omega$ is simply connected, we can write
\beq
\label{gradw}
JU=Dw
\eeq
for a certain smooth function $w$ on $\Omega$. Using that $J$ is an isometry, we have
\beq
\label{eq11}
|Dw|^2=\frac{|Du|^2}{1+|Du|^2}<1.
\eeq
By \rf{eq11}, the function $w$ defines a spacelike graph over $\Omega$, and using that $J(Dw)=-U$ we can also see that $JW=D(-u)$, where
\[
W=\frac{Dw}{\sqrt{1-|Dw|^2}}.
\]
Therefore, $\omega_{JW}$ is also closed on $\Omega$ and, equivalently, $\mathrm{Maximal}[w]=0$. Moreover, the relation $JU=Dw$ jointly with \rf{eq11} implies that $Du$ is parallel if and only if $Dw$ is parallel. A very similar argument, starting now with a non-trivial solution of $\mathrm{Maximal}[w]=0$ with $|Dw|^2<1$ on $\Omega$, produces a non-trivial solution of $\mathrm{Minimal}[u]=0$
\end{proof}

The interest of Theorem\rl{duality} relays on the fact that it allows us to construct new solutions to the maximal
surface equation from known solutions to the minimal surface equation, and viceversa. In particular, as an application of
it we are able to construct counterexamples which show that our Calabi-Bernstein results are no longer true without
the hypothesis $K_M\geq 0$. To see it, let us consider the half-plane model of the hyperbolic plane $\H^2$; that is,
$$
\H^2=\{ x=(x_1,x_2)\in\R^2 : x_2>0 \}
$$
endowed with the complete metric
$$
\g{}{}_{\H^2}=\frac{1}{x_2^2}(dx_1^2+dx_2^2),
$$
conformal to the flat Euclidean metric. For a given smooth function $u=u(x)\in\mathcal{C}^\infty(\H^2)$, we have
that its hyperbolic gradient $Du$ in $\H^2$ and its Euclidean gradient $D_\mathrm{o}u$ in $\R^2$ are related by
\beq
\label{grads}
Du(x)=x_2^2D_\mathrm{o}u(x), \quad x=(x_1,x_2),
\eeq
and then
\beq
\label{eq14}
|Du(x)|^2=x_2^2|D_\mathrm{o}u(x)|^2_\mathrm{o},
\eeq
where $|\cdot|$ and $|\cdot|_\mathrm{o}$ denote, respectively, the norm of a vector field in $\H^2$ and in $\R^2$.
In particular,
\beq
\label{eq12}
\frac{Du(x)}{\sqrt{1+|Du(x)|^2}}=\frac{x_2^2D_\mathrm{o}u(x)}{\sqrt{1+x_2^2|D_\mathrm{o}u(x)|^2_\mathrm{o}}}.
\eeq
The divergence Div of the hyperbolic metric and the divergence $\mathrm{Div}_\mathrm{o}$ of the Euclidean metric
are related by
$$
\mathrm{Div}=\mathrm{Div}_\mathrm{o}-\frac{2}{x_2}dx_2.
$$
By \rf{eq12}, this implies that
\beq
\label{eq13}
\mathrm{Minimal}[u]=\frac{x_2^2\Delta_\mathrm{o}u}{\sqrt{1+x_2^2|D_\mathrm{o}u|^2_\mathrm{o}}}
-\frac{x_2^2}{(1+x_2^2|D_\mathrm{o}u|^2_\mathrm{o})^{3/2}}
\left(x_2u_{x_2}|D_\mathrm{o}u|^2_\mathrm{o}+x_2^2Q(u)\right),
\eeq
where $\Delta_\mathrm{o}$ stands for the Euclidean Laplacian, and
$$
Q(u)=u^2_{x_1}u_{x_1x_1}+2u_{x_1}u_{x_2}u_{x_1x_2}+u^2_{x_2}u_{x_2x_2}.
$$

From here, it is a straightforward computation to check that the function
\beq
\label{minimal1}
u(x_1,x_2)=\log(x_1^2+x_2^2)
\eeq
defines a
non-trivial entire minimal graph in $\H^2\times\R$. Actually, $u$ satisfies
$$
u_{x_2}=\frac{2x_2}{x_1^2+x_2^2}, \qquad \Delta_\mathrm{o}u=0, \qquad Q(u)=\frac{-8}{(x_1^2+x_2^2)^2},
$$
and
\beq
\label{eq15}
|D_\mathrm{o}u|^2_\mathrm{o}=\frac{4}{x_1^2+x_2^2},
\eeq
which by \rf{eq13} gives $\mathrm{Minimal}[u]=0$. Analogously, the function
\beq
\label{minimal2}
u(x_1,x_2)=\frac{x_1}{x_1^2+x_2^2}
\eeq
defines also a non-trivial entire minimal graph in $\H^2\times\R$, because it satisfies
$$
u_{x_2}=\frac{-2x_1x_2}{(x_1^2+x_2^2)^2}, \qquad \Delta_\mathrm{o}u=0, \qquad Q(u)=\frac{2x_1}{(x_1^2+x_2^2)^4},
$$
and
\beq
\label{eq15bis}
|D_\mathrm{o}u|^2_\mathrm{o}=\frac{1}{(x_1^2+x_2^2)^2},
\eeq
which by \rf{eq13} gives $\mathrm{Minimal}[u]=0$.
As far as we know, these examples, which are due to Montaldo and Onnis \cite[Examples 2.3 and 2.4]{MO} (see also
\cite{O}), are the first explicit examples of non-trivial entire minimal graphs in $\H^2\times\R$. Previous existence
results of such minimal graphs have been given in \cite{DH} and \cite{NR}.

\begin{example}
\label{ex1}
From our Theorem\rl{duality} and the entire minimal graph defined by the function \rf{minimal1}, we know that
there exists a smooth function $w\in\mathcal{C}^\infty(\H^2)$ which determines
a non-trivial entire maximal graph in $\H^2\times\R_1$. This shows that the assumption $K_M\geq 0$ in Theorem\rl{CBth2}
and Corollary\rl{CBco2} is necessary. Moreover, we claim that the entire maximal graph determined by $w$ is also complete,
in the sense that the induced metric on $\H^2$ via the graph is complete. To see it, denote by $\g{}{}$ that metric, which
is given by
$$
\g{}{}=\g{}{}_{\H^2}-dw^2
$$
(see equation \rf{gu}). Then, for every tangent vector field $X$ on $\H^2$ we have that
$$
\g{X}{X}=\g{X}{X}_{\H^2}-X(w)^2=\g{X}{X}_{\H^2}-\g{X}{Dw}^2_{\H^2},
$$
and using Cauchy-Schwarz inequality here we get that
\beq
\label{eq16}
\g{X}{X}\geq\g{X}{X}_{\H^2}(1-|Dw|^2).
\eeq
From \rf{eq11} we know that
$$
1-|Dw|^2=\frac{1}{1+|Du|^2}.
$$
In our case, by \rf{eq14} and \rf{eq15} this gives
$$
1-|Dw|^2=\frac{1}{1+x_2^2|D_\mathrm{o}u|^2_\mathrm{o}}=\frac{x_1^2+x_2^2}{x_1^2+5x_2^2}\geq\frac{1}{5},
$$
which jointly with \rf{eq16} yields
$$
\g{}{}\geq\frac{1}{5}\g{}{}_{\H^2}.
$$
As a consequence, the metric $\g{}{}$ is complete on $\H^2$, as claimed. This shows that the assumption $K_M\geq 0$ in
Theorem\rl{CBth} and Corollary\rl{CBco} is also necessary.

It is even possible to get explicitly $w$. In fact, observe that the hyperbolic gradient $Dw$ is given by \rf{gradw}.
Then, from the relation \rf{grads} between the hyperbolic and the Euclidean gradients, this implies that
\[
D_\mathrm{o}w=\left(w_{x_1},w_{x_2}\right)=
\frac{1}{\sqrt{1+x_2^2|D_\mathrm{o}u|^2_\mathrm{o}}}J_\mathrm{o}(D_\mathrm{o}u),
\]
where $J_\mathrm{o}$ denotes the positive $\pi/2$-rotation on the plane, that is,
\[
J_\mathrm{o}(D_\mathrm{o}u)=J_\mathrm{o}(u_{x_1},u_{x_2})=\left(-u_{x_2},u_{x_1}\right)
\]
Therefore, by \rf{eq15bis} we conclude that
\beq
\label{grafo01}
w_{x_1}=-\sqrt{\frac{x_1^2+x_2^2}{x_1^2+5x_2^2}}u_{x_2}=
-\frac{2x_2}{\sqrt{(x_1^2+x_2^2)(x_1^2+5x_2^2)}},
\eeq
and
\beq
\label{grafo02}
w_{x_2}=\sqrt{\frac{x_1^2+x_2^2}{x_1^2+5x_2^2}}u_{x_1}=
\frac{2x_1}{\sqrt{(x_1^2+x_2^2)(x_1^2+5x_2^2)}}.
\eeq
We can explicitly integrate equations \rf{grafo01} and \rf{grafo02} obtaining
\[
w(x_1,x_2)=i\frac{2}{\sqrt{5}}F\left(\arcsin\left(i\frac{x_1}{x_2}\right),\frac{1}{\sqrt{5}}\right)+c,
\]
where $c$ is a real constant, $i$ stands for the imaginary unit and $F(\phi,k)$ stands for the elliptic integral of the
first kind with elliptic modulus $k$ and Jacobi amplitude $\phi$. See Figure\rl{ex1fig} for a picture of this graph in
the case $c=0$.

\begin{figure}[h]
\begin{center}
\includegraphics[width=6cm]{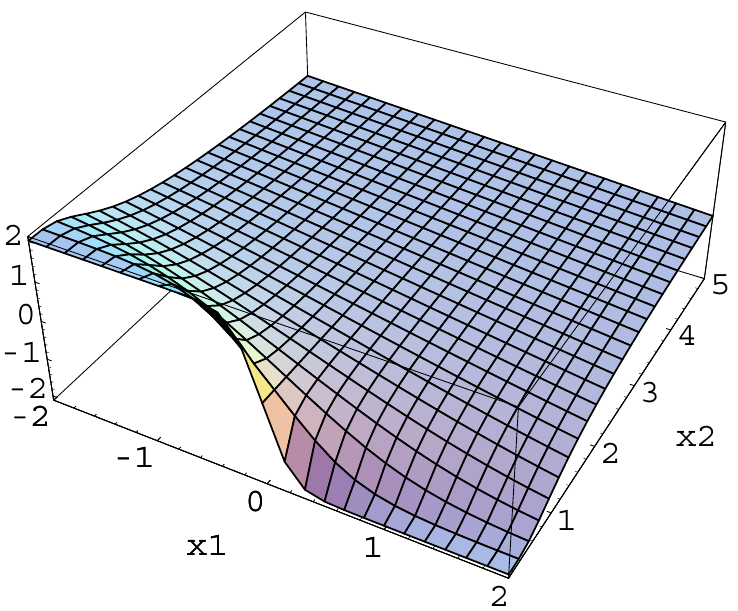}
\caption{\label{ex1fig}The entire complete maximal graph in $\mathbb{H}^2\times\R_1$ given in Example\rl{ex1}.}
\end{center}
\end{figure}

\end{example}

\begin{example}
\label{ex2}
Similarly, the non-trivial entire minimal graph in $\H^2\times\R$ defined by the function \rf{minimal2} gives rise via
Theorem\rl{duality} to
another non-trivial entire maximal graph in the Lorentzian product $\H^2\times\R_1$. In contrast to Example\rl{ex1}, this example is not
complete. To see it, let $w\in\mathcal{C}^\infty(\H^2)$ stand for the smooth function defining this entire maximal graph,
which we denote by $\Sigma(w)$.
In an analogous way as in Example\rl{ex1} we can compute
\beq
\label{grafo1}
w_{x_1}=-\frac{(x_1^2+x_2^2)}{\sqrt{(x_1^2+x_2^2)^2+x_2^2}}u_{x_2}=
\frac{2x_1x_2}{(x_1^2+x_2^2)\sqrt{(x_1^2+x_2^2)^2+x_2^2}},
\eeq
and
\beq
\label{grafo2}
w_{x_2}=\frac{(x_1^2+x_2^2)}{\sqrt{(x_1^2+x_2^2)^2+x_2^2}}u_{x_1}=
\frac{x_2^2-x_1^2}{(x_1^2+x_2^2)\sqrt{(x_1^2+x_2^2)^2+x_2^2}}.
\eeq
Let $\alpha:(0,1)\rightarrow\Sigma(w)$ be the divergent curve in $\Sigma(w)$ given by
\[
\alpha(s)=(0,s,w(0,s)).
\]
Then $\alpha'(s)=(0,1,w_{x_2}(0,s))$ and
\[
\|\alpha'(s)\|^2=\frac{1}{s^2}-w_{x_2}(0,s)^2=\frac{1}{1+s^2},
\]
which implies that $\alpha$ has finite length, because of
$$
\int_0^1\|\alpha'(s)\|ds=\int_0^1\frac{ds}{\sqrt{1+s^2}}=\mathrm{arcsinh}(1)=\log{(1+\sqrt{2})}.
$$
As a consequence, $\Sigma(w)$ is not complete. This fact is particularly interesting. Let us recall that such
circumstance cannot occur in the
Lorentz-Minkowski space $\mathbb{R}^3_1$, since by a result of Cheng and Yau \cite{CY}, surfaces with constant mean curvature and
closed in $\mathbb{R}^3_1$ are necessarily complete.

We can explicitly integrate equations \rf{grafo1} and \rf{grafo2}, getting
\[
w(x_1,x_2)=\ln\left(\frac{x_1^2+x_2^2}{2(x_2+\sqrt{x_2^2+(x_1^2+x_2^2)^2})}\right)+c,
\]
where $c$ is a real constant. Below we exhibit a picture of this graph for $c=0$.
\begin{figure}[h]
\begin{center}
\includegraphics[width=6cm]{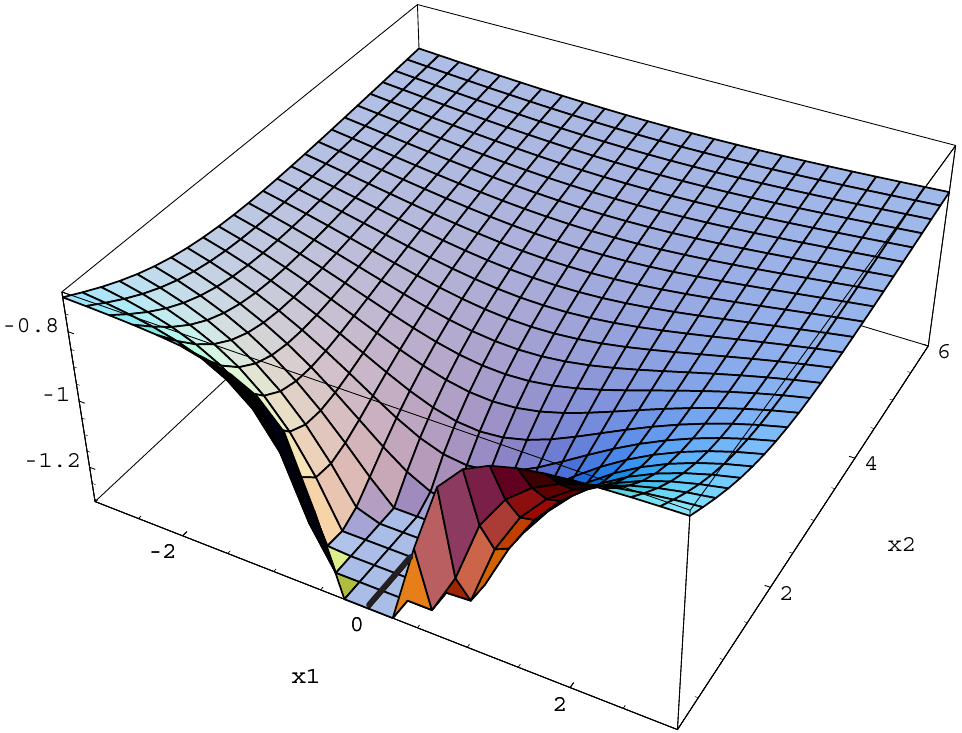}
\caption{The entire non-complete maximal graph in $\mathbb{H}^2\times\R_1$ given in Example\rl{ex2}.}
\end{center}
\end{figure}

\end{example}

See also \cite{Al} for further examples of complete and non-complete entire maximal graphs in $\H^2\times\R_1$
obtained by the first author by looking for explicit solutions of the partial differential equation \rf{zmc} on
$\H^2$.

\section*{Acknowledgements}
This work was written while the second author was visiting the Institut des Hautes \'{E}tudes Scientifiques (IH\'{E}S) in
Bures-sur-Yvette, France. He thanks IH\'{E}S for its hospitality and support. The authors would like to heartily thank M. Dajczer, E. Garc\'\i a-R\'\i o, B. Palmer, M.
S\'{a}nchez and R. Souam for several suggestions and useful comments during the preparation of this paper.
They also thank to the referee for valuable suggestions which improved the paper. This paper is part of the PhD thesis
of A.L. Albujer \cite{Al2}, which was defended in Universidad de Murcia, Spain, in November 2008.

\bibliographystyle{amsplain}

\begin{thebibliography}{10}

\bibitem{Ah} L.V. Ahlfors, \textit{Sur le type d'une surface de Riemann},
C.R. Acad. Sc. Paris {\bf 201} (1935), 30--32.

\bibitem{Al} A.L. Albujer,
\textit{New examples of entire maximal graphs in $\H^2\times\R_1$}, Differential Geom. Appl.
\textbf{26} (2008), 456--462.

\bibitem{Al2} A.L. Albujer,
\textit{Geometr\'\i a global de superficies espaciales en espacios producto lorentzianos}, Ph.D. thesis, Universidad de Murcia, Spain
(2008).

\bibitem{AP1} L.J. Al\'\i as and B. Palmer,
\textit{A duality result between the minimal surface equation and the maximal surface
equation}, An. Acad. Bras. Ci\^{e}nc. \textbf{73} (2001), 161--164.

\bibitem{AP2} L.J. Al\'\i as and B. Palmer,
\textit{On the Gaussian curvature of maximal surfaces and the Calabi-Bernstein theorem},
Bull. London Math. Soc. {\bf 33} (2001), 454--458.

\bibitem{ARS1} L.J. Al\'\i as, A. Romero and M. S\'{a}nchez,
\textit{Uniqueness of complete spacelike hypersurfaces of constant mean curvature in generalized Robertson-Walker
spacetimes}, Gen. Relativity Gravitation \textbf{27} (1995), 71--84.

\bibitem{BF} D. Brill and F. Flaherty,
\textit{Isolated maximal surfaces in spacetime}, Comm. Math. Phys.
\textbf{50} (1976), 157--165.


\bibitem{Ca} E. Calabi, \textit{Examples of Bernstein problems for some nonlinear equations}.
1970  Global Analysis (Proc. Sympos. Pure Math., Vol. XV, Berkeley, Calif., 1968)  pp. 223--230 Amer. Math. Soc., Providence, R.I.

\bibitem{CY} S.Y. Cheng and S.T. Yau,
\textit{Maximal space-like hypersurfaces in the Lorentz-Minkowski spaces},
Ann. of Math. (2) {\bf 104} (1976), 407--419.

\bibitem{Ch} S.S. Chern, \textit{Simple proofs of two theorems in minimal surfaces},
Enseign. Math. II. S\'{e}r. {\bf 15} (1969), 53--61.

\bibitem{DH} D.M. Duc and N.V. Hieu, \textit{Graphs with prescribed mean
curvature on Poincar\'{e} disk}, Bull. London Math. Soc. {\bf 27} (1995), 353--358.

\bibitem{ER1} F.J.M. Estudillo and A. Romero,
\textit{On maximal surfaces in the $n$-dimensional Lorentz-Minkowski space},
Geom. Dedicata {\bf 38}  (1991), 167--174.

\bibitem{ER2} F.J.M. Estudillo and A. Romero,
\textit{Generalized maximal surfaces in Lorentz-Minkowski space $L\sp 3$},
Math. Proc. Cambridge Philos. Soc. {\bf 111} (1992), 515--524.

\bibitem{ER3} F.J.M. Estudillo and A. Romero,
\textit{On the Gauss curvature of maximal surfaces in the $3$-dimensional Lorentz-Minkowski space},
Comment. Math. Helv. {\bf 69} (1994),  1--4.

\bibitem{FM} I. Fern\'{a}ndez and P. Mira,
\textit{Complete maximal surfaces in static Robertson-Walker $3$-spaces},
Gen. Relativ. Gravit. {\bf 39} (2007),  2073--2077.


\bibitem{Fr} T. Frankel,
\textit{Applications of Duschek's formula to cosmology and minimal surfaces}, Bull. Amer. Math. Soc.
\textbf{81} (1975), 579--582.

\bibitem{Hu} A. Huber, \textit{On subharmonic functions and differential
geometry in the large}, Comment. Math. Helv. {\bf 32} (1957), 13--72.

\bibitem{Ko} O. Kobayashi,
\textit{Maximal surfaces in the $3$-dimensional Minkowski space $L\sp{3}$},
Tokyo J. Math. {\bf 6} (1983), 297--309.

\bibitem{KN} S. Kobayashi and K. Nomizu, Foundations of Differential Geometry, Vol. II,
Interscience, New York, 1969.

\bibitem{LR1} J.M. Latorre and A. Romero,
\textit{New examples of Calabi-Bernstein problems for some nonlinear equations}, Differential Geom. Appl.
\textbf{15} (2001), 153--163.

\bibitem{MT} J.E. Marsden and F.J. Tipler,
\textit{Maximal hypersurfaces and foliations of constant mean curvature in general relativity}, Phys. Rep.
\textbf{66} (1980), 109--139.


\bibitem{Mc} L.V. McNertey, \textit{One-parameter families of surfaces with constant curvature in Lorentz 3-space}, Ph.D. thesis, Brown University (USA), 1980.

\bibitem{MO} S. Montaldo and I.I. Onnis, \textit{A note on surfaces in $\H^2\times\R$},
Boll. Unione Mat. Ital. Sez. B Artic. Ric. Mat. (8) {\bf 10} (2007), 939--950.

\bibitem{NR} B. Nelli and H. Rosenberg, \textit{Minimal surfaces in $\H^2\times\R$},
Bull. Braz. Math. Soc. (N.S.) {\bf 33} (2002), 263--292.

\bibitem{O} I.I. Onnis, \textit{Superficies em certos espacos homogeneos tridimensionais}, Ph.D. thesis, Universidade
Estadual de Campinas (Brazil), 2005.

\bibitem{Ro} A. Romero,
\textit{Simple proof of Calabi-Bernstein's theorem on maximal surfaces}, Proc. Amer. Math. Soc. \textbf{124} (1996),
1315--1317.

\end{thebibliography}

\end{document}